# Some Bounds on the Lehmer Index

Xiaodong Zhuang
Qingdao University, Shandong Province, China
xdzhuang@qdu.edu.cn

## 1. Introduction

A Lehmer number is a positive composite integer satisfying $\varphi(n)|n-1$[1]. The Lehmer conjecture on Euler totient states that there is no Lehmer number[1]. Currently, this conjecture is an unsolved problem. Existing work is mainly focused on raising the lower bounds of the size and factor number of Lehmer number, as well as other constraints that a Lehmer number should satisfy[2,3,4]. These constraints progressively narrow the range of Lehmer number.

In this paper, the relationship is studied between the Lehmer index and the reciprocal sum of the prime factors of a Lehmer number. Some upper bounds and a lower bound are proposed for the Lehmer index. Another constraint on Lehmer index is also proposed if a Lehmer number has 3 as its prime factor.

## 2. Some upper bounds of the Lehmer index

Lehmer has proven that a Lehmer number must be square-free and odd, and has at least 7 distinct prime factors[1]. Assume n is a Lehmer number, and $n = p_1 p_2 \dots p_k$. Its Euler function is $\varphi(n) = (p_1-1)(p_2-1)\dots(p_k-1)$. Then n's factors satisfy:

$$\left(\prod_{i=1}^{k} p_i\right) - 1 = u \cdot \prod_{i=1}^{k}(p_i - 1) \tag{1}$$

where $p_1, p_2, \dots p_k$ are n's prime factors, and u is named the Lehmer index of n[4]. In the following, the upper bound of u is derived from Equation (1).

First, divide both sides of (1) by nu and rearrange to obtain:

$$\frac{1}{u}\left(1-\frac{1}{n}\right) = \prod_{i=1}^{k}\left(1-\frac{1}{p_i}\right)$$

Then take the natural logarithm of both sides of the above, we have:

$$-ln(u) + \ln\left(1-\frac{1}{n}\right) = \sum_{i=1}^{k} ln\left(1-\frac{1}{p_i}\right) \tag{2}$$

Using the Taylor series expansion of ln(1-x):

$$ln(1-x) = -x - \frac{x^2}{2} - \frac{x^3}{3} - \cdots = -\sum_{j=1}^{+\infty}\frac{x^j}{j} \qquad x \in [0,1) \tag{3}$$

we can turn the right-hand side of Equation (2) into:

$$\sum_{i=1}^{k} ln\left(1-\frac{1}{p_i}\right) = -\sum_{i=1}^{k}\sum_{j=1}^{+\infty}\left(\frac{1}{j}\cdot\frac{1}{{p_i}^j}\right) = -\sum_{j=1}^{+\infty}\left(\frac{1}{j}\sum_{i=1}^{k}\frac{1}{{p_i}^j}\right) \tag{4}$$

Let $s_j$ denote the reciprocal sum of the j-th power of $p_i$:

$$s_j = \sum_{i=1}^{k}\frac{1}{{p_i}^j}$$

Substitute $s_j$ into (4), we have:

$$-ln(u) + \ln\left(1-\frac{1}{n}\right) = -\sum_{j=1}^{+\infty}\frac{s_j}{j}$$

and

$$ln(u) - \ln\left(1-\frac{1}{n}\right) = \sum_{j=1}^{+\infty}\frac{s_j}{j} \tag{5}$$

Because $\ln(1-1/n) < 0$, we have:

$$ln(u) = \sum_{j=1}^{+\infty}\frac{s_j}{j} + \ln\left(1-\frac{1}{n}\right) < \sum_{j=1}^{+\infty}\frac{s_j}{j} \tag{6}$$

The term $\ln(1-1/n)$ which is dropped in the above inequality estimate is very small in magnitude. Pinch proved that the Lehmer number is not less than $10^{30}$. Hence $|\ln(1-1/n)| \approx 1/n \leq 10^{-30}$.

From (6), an upper bound of Lehmer index u can be derived. Let $p_1$ denote the smallest prime factor of n. Then $p_i > p_1$, $1/p_i < 1/p_1$ for $1 < i \leq k$. For $s_j$ $(j \geq 2)$:

$$s_j = \sum_{i=1}^{k}\frac{1}{{p_i}^j} < \sum_{i=1}^{k}\frac{1}{p_i\cdot{p_1}^{j-1}} = \frac{1}{{p_1}^{j-1}}\sum_{i=1}^{k}\frac{1}{p_i} = \frac{s_1}{{p_1}^{j-1}} \tag{7}$$

where $s_1$ is the reciprocal sum of n's prime factors. Substitute the above into (6), we have:

$$ln(u) < \sum_{j=1}^{+\infty}\frac{s_j}{j} < \sum_{j=1}^{+\infty}\frac{s_1}{j\cdot{p_1}^{j-1}} = s_1\cdot p_1\cdot\sum_{j=1}^{+\infty}\frac{1}{j\cdot{p_1}^j} \tag{8}$$

According to the Taylor series of (3), we have:

$$-ln\left(1-\frac{1}{p_1}\right) = \sum_{j=1}^{+\infty}\frac{1}{j\cdot{p_1}^j}$$

Substitute the above equation into (8), we have:

$$ln(u) < -s_1\cdot p_1\cdot ln\left(1-\frac{1}{p_1}\right)$$

Rearrange the above to obtain:

$$ln(u) < s_1\cdot ln\left(\left(1+\frac{1}{p_1-1}\right)^{p_1}\right) \tag{9}$$

Exponentiating both sides of the above, an upper bound of u is obtained:

$$u < e^{ln\left(\left(1+\frac{1}{p_1-1}\right)^{p_1}\right)\cdot s_1} = \left(\left(1+\frac{1}{p_1-1}\right)^{p_1}\right)^{s_1} \quad (10)$$

Because $p_1$ is an odd prime, $p_1 \geq 3$, and $1/p_1 \leq 1/3$. Hence, a more general upper bound of $s_1$ can be obtained from (9) and (10):

$$u < \left(\frac{27}{8}\right)^{s_1} \quad (11)$$

It is obvious that (10) is a tighter upper bound than (11). Moreover, since $p_i \geq 3$ there is natural relationship between $s_1$ and k (k is the prime factor number of n):

$$s_1 = \sum_{i=1}^{k}\frac{1}{p_i} < \sum_{i=1}^{k}\frac{1}{3} = \frac{k}{3}$$

Hence another upper bound relating $s_1$ to k is:

$$u < \left(\frac{27}{8}\right)^{s_1} < \left(\frac{27}{8}\right)^{\frac{k}{3}} = 1.5^k \quad (12)$$

Obviously, (12) is a weaker upper bound than (11).

## 3. A lower bound of the Lehmer index

Separate the item of $s_1$ from the sum on the right-hand side in (5), we have:

$$ln(u) - \ln\left(1-\frac{1}{n}\right) = s_1 + \sum_{j=2}^{+\infty}\frac{s_j}{j}$$

Rearrange the above, we have:

$$ln(u) = s_1 + \sum_{j=2}^{+\infty}\frac{s_j}{j} + \ln\left(1-\frac{1}{n}\right) > s_1 + \ln\left(1-\frac{1}{n}\right)$$

Exponentiating the left and right-hand sides of the above, we have:

$$u > e^{s_1+\ln\left(1-\frac{1}{n}\right)} = e^{s_1}\left(1-\frac{1}{n}\right)$$

Using Pinch's result: $1/n \leq 10^{-30}$ [ ], a lower bound of u is obtained:

$$u > e^{s_1}\cdot(1-10^{-30}) \quad (13)$$

## 4. The constraint on Lehmer index brought by the prime factor 3

If 3 is one prime factor of n, it brings extra constraint on the Lehmer index. Substitute $p_1 = 3$ into Equation (1), we have:

$$3\cdot\left(\prod_{i=2}^{k}p_i\right) - 1 = u\cdot(3-1)\cdot\prod_{i=2}^{k}(p_i-1)$$

and

$$3 \cdot \left( \prod_{i=2}^{k} p_i \right) = 2u \cdot \prod_{i=2}^{k} (p_i - 1) + 1 \tag{14}$$

Because n and (n-1) are coprime, any factor of n is not the factor of (n-1). Then Equation (1) indicates that any factor of n is not the factor of u and $(p_i - 1)$, $i = 1, 2, \ldots, k$. Since $p_1 = 3$, $(p_i - 1)$ is not divisible by 3 (i.e. $p_i - 1 \not\equiv 0 \pmod 3$). On the other hand, the prime factors of n are distinct, hence $p_i$ ($i = 2, \ldots, k$) is not divisible by 3. That is to say, $p_i - 1 \not\equiv 2 \pmod 3$, $i = 2, \ldots, k$. Therefore, we have: $p_i - 1 \equiv 1 \pmod 3$. Then let's analyze both sides of Equation (14) modulo 3. The left-hand side is divisible by 3, which forces $u \equiv 1 \pmod 3$ in order that the right-hand side of (14) can also be divisible by 3. Therefore, the following theorem holds:

Theorem If 3 is a prime factor of a Lehmer number, its Lehmer index must be congruent to 1 modulo 3, thus the smallest value of its Lehmer index is 4.

## 5. Discussion

The upper and lower bounds proposed for Lehmer index provide more constraints on Lehmer number, which can be combined with other methods to derive more results. The combination with number theory analysis is a direction in future work.